\pdfoutput=1
\RequirePackage{ifpdf}
\ifpdf
\documentclass[pdftex]{sigma}
\else
\documentclass{sigma}
\fi

\newcommand{\pos}{{\mathrm{pos}}}

\newcommand{\cc}{cc}
\newcommand{\bsi}{{\overline{\sigma}}}
\newcommand{\bta}{{\overline{\tau}}}

\newcommand{\C}{\mathbb{C}}

\newcommand{\lam}{\boldsymbol{\lambda}}
\newcommand{\num}{\boldsymbol{\nu}}
\newcommand{\mum}{\boldsymbol{\mu}}
\newcommand{\alm}{\boldsymbol{\alpha}}
\newcommand{\bem}{\boldsymbol{\beta}}
\newcommand{\gam}{\boldsymbol{\gamma}}
\newcommand{\dem}{\boldsymbol{\delta}}
\newcommand{\tem}{\boldsymbol{\theta}}

\newtheorem*{theo-intro}{Theorem}
\newtheorem{Theorem}{Theorem}
\newtheorem{Lemma}[Theorem]{Lemma}
\newtheorem{Proposition}[Theorem]{Proposition}

{\theoremstyle{definition}
\newtheorem*{Remarks}{Remarks}
\newtheorem*{Example}{Example}
}

\begin{document}

\allowdisplaybreaks

\renewcommand{\thefootnote}{$\star$}

\renewcommand{\PaperNumber}{039}

\FirstPageHeading

\ShortArticleName{Fusion Procedure for Cyclotomic Hecke Algebras}

\ArticleName{Fusion Procedure for Cyclotomic Hecke Algebras\footnote{This paper is a~contribution to the Special Issue in honor of Anatol
Kirillov and Tetsuji Miwa.
The full collection is available at \href{http://www.emis.de/journals/SIGMA/InfiniteAnalysis2013.html}
{http://www.emis.de/journals/SIGMA/InfiniteAnalysis2013.html}}}

\Author{Oleg V.~OGIEVETSKY~$^{\dag^1\dag^2\dag^3}$ and Lo\"{\i}c POULAIN D'ANDECY~$^{\dag^4}$}

\AuthorNameForHeading{O.V.~Ogievetsky and L.~Poulain d'Andecy}

\Address{$^{\dag^1}$~Center of Theoretical Physics,
Aix Marseille Universit\'e, CNRS,\\
\hphantom{$^{\dag^1}$}~UMR 7332, 13288 Marseille, France}
\EmailDD{\href{mailto:oleg@cpt.univ-mrs.fr}{oleg@cpt.univ-mrs.fr}}
\Address{$^{\dag^2}$~Universit\'e de Toulon, CNRS, UMR 7332, 83957 La Garde, France}
\Address{$^{\dag^3}$~On leave of absence from P.N.~Lebedev Physical Institute,\\
\hphantom{$^{\dag^3}$}~Leninsky Pr.~53, 117924 Moscow, Russia}

\Address{$^{\dag^4}$~Mathematics Laboratory of Versailles, LMV, CNRS UMR 8100, \\
\hphantom{$^{\dag^4}$}~Versailles Saint-Quentin University, 45 avenue des Etas-Unis, \\
\hphantom{$^{\dag^4}$}~78035 Versailles Cedex, France}
\EmailDD{\href{mailto:L.B.PoulainDAndecy@uva.nl}{L.B.PoulainDAndecy@uva.nl}}

\ArticleDates{Received September 28, 2013, in f\/inal form March 29, 2014; Published online April 01, 2014}

\Abstract{A complete system of primitive pairwise orthogonal idempotents for cyclotomic Hecke algebras is constructed by consecutive
evaluations of a rational function in several variables on quantum contents of multi-tableaux. This function is a product of two terms, one of which depends only on the shape of the multi-tableau and is proportional to the inverse of the corresponding Schur element.}


\Classification{20C08; 05E10}

\renewcommand{\thefootnote}{\arabic{footnote}}
\setcounter{footnote}{0}

\section{Introduction}

This article is a continuation of the article \cite{OPdA} on the fusion procedure for the complex ref\/lection groups $G(m,1,n)$.
The cyclotomic Hecke algebra $H(m,1,n)$, introduced in \cite{AK,Bro-M,C2}, is a natural f\/lat
deformation of the group ring of the complex ref\/lection group $G(m,1,n)$.

In \cite{OPdA}, a fusion procedure, in the spirit of \cite{Mo}, for the complex ref\/lection groups $G(m,1,n)$ is suggested:
a complete system of primitive pairwise orthogonal idempotents for the groups $G(m,1,n)$ is obtained by consecutive evaluations of a rational function in several variables with values in the group ring $\C G(m,1,n)$. This approach to the fusion procedure relies on the existence of a maximal commutative set of elements of $\C G(m,1,n)$ formed by the Jucys--Murphy elements.

Jucys--Murphy elements for the cyclotomic Hecke algebra $H(m,1,n)$ were introduced in \cite{AK} and were used in \cite{OPdA1} to develop an inductive approach to the representation theory of the chain of the algebras $H(m,1,n)$. In the generic setting
or under certain restrictions on the parameters of the algebra $H(m,1,n)$
(see Section \ref{sec-def} for precise def\/initions), the Jucys--Murphy elements form a maximal commutative set in the algebra $H(m,1,n)$.

A complete system of primitive pairwise orthogonal idempotents of the algebra $H(m,1,n)$ is indexed by the set of standard $m$-tableaux of size $n$. We formulate here the main result of the article.
Let $\lam$ be an $m$-partition of size $n$ and ${\cal{T}}$ be a standard $m$-tableau of shape $\lam$.

\begin{theo-intro}\label{prop-fus1-intro}
The idempotent $E_{{\cal{T}}}$ of $H(m,1,n)$
corresponding to the standard $m$-tableau ${\cal{T}}$ of shape $\lam$ can be obtained by the following consecutive evaluations
\begin{gather}\label{eq-idem-fin1-intro}
E_{{\cal{T}}}=\mathsf{F}_{\lam}\Phi(u_1,\dots,u_{n})\Bigr\rvert_{u_1=c_1}\cdots\Bigr\rvert_{u_{n-1}=c_{n-1}}\Bigr\rvert_{u_{n}=c_{n}}\ .
\end{gather}
\end{theo-intro}

Here $\Phi(u_1,\dots,u_{n})$ is a rational function with values in the algebra $H(m,1,n)$, $\mathsf{F}_{\lam}$ is an element of the base ring and $c_1,\dots,c_{n}$ are the quantum contents of the $m$-nodes of ${\cal{T}}$.

The classical limit of our fusion procedure for algebras $H(m,1,n)$
reproduces the fusion procedure of \cite{OPdA} for the complex ref\/lection groups $G(m,1,n)$.
 For  $\C G(m,1,n)$, the variables of the rational function are split into
two parts, one is related to the position of the $m$-node (its place in the $m$-tuple) and the other one~--
to the classical content of the $m$-node.
The position variables can be evaluated simultaneously while the classical content variables have then to be
evaluated consequently from 1 to $n$. For the algebra $H(m,1,n)$, the information about positions and classical contents is fully contained in the quantum contents, and now the function $\Phi$ depends on only one set of variables.

Remarkably, the coef\/f\/icient $\mathsf{F}_{\lam}$ appearing in (\ref{eq-idem-fin1-intro})
depends only on the shape $\lam$ of the standard $m$-tableau ${\cal{T}}$ (cf.\ with the more delicate
fusion  procedure for the Birman--Murakami--Wenzl algebra \cite{IMO3}). In the classical limit, this coef\/f\/icient depends only on the usual hook length,
see~\cite{OPdA}. However, in the deformed situation, the calculation of $\mathsf{F}_{\lam}$ needs a non-trivial generalization of the hook length. It appears that the coef\/f\/icient $\mathsf{F}_{\lam}$ is proportional to the inverse of the {\it Schur element} (corresponding to the $m$-partition $\lam$) associated to a specif\/ic
symmetrizing form on the algebra~$H(m,1,n)$
(see~\cite{GIM,M} for a calculation of these Schur elements and~\cite{CJ} for an expression in terms of generalized hook lengths);
for more precise statements, we refer to~\cite{OPdA2} where we calculate, using the fusion formula presented here,
weights of certain central forms and in particular of these Schur elements.

For $m=1$, the cyclotomic Hecke algebra $H(1,1,n)$ is the Hecke algebra of type A and
our fusion procedure reduces to the fusion procedure for the Hecke algebra in~\cite{IMO}.
The factors in the rational function are arranged in \cite{IMO} in such a way that there is
a product of ``Baxterized'' generators on one side and a product of non-Baxterized generators on the other side.
{}For $m>1$ a rearrangement, as for the type~A, of the rational function appearing in~(\ref{eq-idem-fin1-intro}) is no more possible.

The additional, with respect to $H(1,1,n)$,  generator of $H(m,1,n)$ satisf\/ies the ref\/lection equation whose ``Baxterization'' is known
\cite{IO2}. But~-- and this is maybe surprising~-- the full Baxterized form is not used in the construction of the rational function in  (\ref{eq-idem-fin1-intro}).
The rational expression involving the additional generator satisf\/ies only a certain limit of the ref\/lection equation with spectral parameters.

The Hecke algebra of type A is the natural quotient of the Birman--Murakami--Wenzl algebra. The fusion procedure, developed in~\cite{IMO3}, for the Birman--Murakami--Wenzl algebra provides a one-parameter family of fusion procedures for the Hecke algebra of type~A. We think that for $m>1$ the
fusion procedure (\ref{eq-idem-fin1-intro}) can be included into a one-parameter family as well.

\section{Def\/initions}\label{sec-def}

\subsection{Cyclotomic Hecke algebra and Baxterized elements}\label{subsec-def}

Let $m\in\mathbb{Z}_{>0}$ and $n\in\mathbb{Z}_{\geq 0}$. Let $q,v_1,\dots,v_m$ be complex numbers with $q\neq 0$.
The cyclotomic Hecke algebra $H(m,1,n+1)$ is the unital associative algebra over $\mathbb{C}$ generated by $\tau$, $\sigma_1,\dots,\sigma_n$ with the def\/ining relations
\begin{alignat*}{3}
& \sigma_i\sigma_{i+1}\sigma_i=\sigma_{i+1}\sigma_i\sigma_{i+1}\qquad  && \text{for $i=1\dots,n-1$},& \nonumber\\
& \sigma_i\sigma_j=\sigma_j\sigma_i \qquad  && \text{for $i,j=1,\dots,n$ such that $|i-j|>1$},& \nonumber\\
& \tau\sigma_1\tau\sigma_1=\sigma_1\tau\sigma_1\tau ,\qquad &&& \nonumber\\
& \tau\sigma_i=\sigma_i\tau \qquad  && \text{for $i>1$},& \nonumber\\
& \sigma_i^2=\big(q-q^{-1}\big)\sigma_i+1 \qquad  && \text{for $i=1,\dots,n$},& \nonumber\\
& (\tau-v_1)\cdots(\tau-v_m)=0  .\qquad &&& 
 \end{alignat*}
We def\/ine $H(m,1,0):=\mathbb{C}$. The cyclotomic Hecke algebras $H(m,1,n)$ form a chain (with respect to~$n$)
of algebras def\/ined by inclusions
$H(m,1,n)\ni \tau,\sigma_1,\dots,\sigma_{n-1}\mapsto  \tau,\sigma_1,\dots,\sigma_{n-1}\in H(m,1,n+1)$ for any $n\geq 0$. These inclusions allow to consider (as it will often be done in the article) elements of $H(m,1,n)$ as elements of $H(m,1,n+n')$ for any $n'=0,1,2,\dots$.

In the sequel we assume the following restrictions on the parameters $q,v_1,\dots,v_m$:
\begin{gather}\label{sesi1}
 1+q^2+\dots+q^{2N}\neq 0\ \ {\textrm{for}}\ N\ \text{such that}\ N\leq n
 , \\
\label{sesi2}
 q^{2i}v_j-v_k\neq 0\ \ {\textrm{for}}\ \ i,j,k \ \ {\textrm{such that}}\ \  j\neq k \ {\textrm{and}}\ -n\leq i \leq n
 , \\
\label{sesi3}
 v_j\neq 0\ \ \text{for}\ \ j=1,\dots,m  .
\end{gather}
The restrictions (\ref{sesi1}), (\ref{sesi2}) are necessary and suf\/f\/icient for the semi-simplicity of
the algebra $H(m,1,n+1)$~\cite[main theorem]{A}.
The restriction~(\ref{sesi3}) is necessary for the maximality of the commutative
set of the Jucys--Murphy elements (as def\/ined in Section~\ref{sec-idem}) \cite[Proposition~3.2]{A}.

Def\/ine the following rational functions in variables $a,b$ with values in $H(m,1,n+1)$:
\begin{gather}\label{bax-s}
\bsi_i(a,b):=\sigma_i+(q-q^{-1})\frac{b}{a-b}  ,\qquad i=1,\dots,n  .
\end{gather}
The functions $\bsi_i$ are called {\it Baxterized} elements and the variables $a$ and $b$ are called {\it spectral parameters}.
These Baxterized elements satisfy the Yang--Baxter equation with spectral para\-me\-ters
\begin{gather*}
\bsi_i(a,b)\bsi_{i+1}(a,c)\bsi_i(b,c)=\bsi_{i+1}(b,c)\bsi_i(a,c)\bsi_{i+1}(a,b) .
\end{gather*}
The following formula will be used later
\begin{gather}\label{baxt-inv1}
 \bsi_i(a,b)\bsi_i(b,a)=\frac{(a-q^2b)(a-q^{-2}b)}{(a-b)^2} \qquad\textrm{for \ \ $i=1,\dots,n$}.
\end{gather}

Let $\mathfrak{p}_i$, $i=1,\dots,m$,
be the eigen-idempotents of $\tau$, $\mathfrak{p}_i:=\prod\limits_{j:j\neq i} (\tau-v_j)/(v_i-v_j)$, so that $\tau\mathfrak{p}_i=v_i\mathfrak{p}_i$, $\mathfrak{p}_i\mathfrak{p}_j=\delta_{ij}\mathfrak{p}_i$, $\sum_i \mathfrak{p}_i=1$ and $\tau=\sum_i v_i\mathfrak{p}_i$.
Let $r$ be an indeterminate.
The resolvent $(r-\tau)^{-1}:=\sum_i(r-v_i)^{-1}\mathfrak{p}_i$ of $\tau$ is an element of
$\mathbb{C}(r)\otimes_{\mathbb{C}} H(m,1,n+1)$.
Def\/ine a rational function $\bta$ with values in $H(m,1,n+1)$:
\begin{gather}\label{bax-t2}
\bta(r):=\frac{(r-v_1)(r-v_2)\cdots(r-v_m)}{r-\tau}=\sum_i\left(\prod_{j:j\neq i} (r-v_j)\!\right)\!\mathfrak{p}_i\in\mathbb{C}[r]\otimes_{\mathbb{C}} H(m,1,n+1)  .\!\!\!
\end{gather}

\begin{Remarks} \textbf{(i)}
The function $\bta(r)$ can be expressed in terms of the complex numbers $a_0,a_1,{\dots},a_m$ def\/ined by
\[
(X-v_1)(X-v_2)\cdots(X-v_m)=a_0+a_1X+\cdots+a_{m}X^{m}  ,
\]
where $X$ is an indeterminate.
Let $\mathfrak{a}_i(r)$, $i=0,\dots,m$, be the polynomials in $r$ given by
\begin{gather}\label{polalicl}
\mathfrak{a}_i(r)=a_i+r a_{i+1}+\cdots+r^{m-i}a_m  \qquad\textrm{for \ \ $i=0,\dots,m$.}
\end{gather}
Using that $r \mathfrak{a}_{i+1}(r)=\mathfrak{a}_{i}(r)-a_i$, for $i=0,\dots,m-1$, it is straightforward to verify that
\begin{gather}
\label{bax-t-inter}
(r-\tau)\sum_{i=0}^{m-1}\mathfrak{a}_{i+1}(r)\tau^i=\mathfrak{a}_0(r)=(r-v_1)(r-v_2)\cdots(r-v_m) .
\end{gather}
It follows from (\ref{bax-t-inter}) that
\begin{gather}\label{bax-t}
\bta(r)=\mathfrak{a}_1(r)+\mathfrak{a}_2(r)\tau+\cdots+\mathfrak{a}_{m}(r)\tau^{m-1}=\sum_{i=0}^{m-1}\mathfrak{a}_{i+1}(r)\tau^i  ,
\end{gather}
For example, for $m=1$, we have $\bta(r)=1$; for $m=2$, we have $\bta(r)=\tau+r-v_1-v_2$; for $m=3$, we have $\bta(r)=\tau^2+(r-v_1-v_2-v_3)\tau+r^2-r(v_1+v_2+v_3)+v_1v_2+v_1v_3+v_2v_3 $.

\textbf{(ii)}
The functions $\bta$ and $\bsi_1$ satisfy the following equation
\begin{gather}\label{ref-baxt}
\bsi_1(a,b)\bta(a)\sigma_1^{-1}\bta(b)=\bta(b)\sigma_1^{-1}\bta(a)\bsi_1(a,b)  .
\end{gather}
Indeed, due to (\ref{baxt-inv1}) and (\ref{bax-t2}), the equality (\ref{ref-baxt}) is equivalent to
\[
(\tau-b)\sigma_1(\tau-a)\bsi_1(b,a)=\bsi_1(b,a)(\tau-a)\sigma_1(\tau-b)  ,
\]
which is proved by a straightforward calculation.
The equation (\ref{ref-baxt}) is a certain (we leave the details to the reader) limit of the usual ref\/lection equation with spectral parameters (see, for example,~\cite{K}).
\end{Remarks}

\subsection[$m$-partitions, m-tableaux and generalized hook length]{$\boldsymbol{m}$-partitions, m-tableaux and generalized hook length}

Let $\lambda\vdash n+1$ be a partition of size $n+1$, that is, $\lambda=(\lambda_1,\dots,\lambda_l)$, where $\lambda _j$, $j=1,\dots,l$,
are positive integers, $\lambda_1\geqslant\lambda_2\geqslant\dots\geqslant\lambda_l$ and $n+1=\lambda_1+\dots+\lambda_l$.
We identify partitions with their Young diagrams: the Young diagram of $\lambda$ is a left-justif\/ied array of rows of
nodes containing $\lambda_j$ nodes in the $j$-th row, $j=1,\dots,l$; the rows are numbered from top to bottom.
For a node $\alpha$ in line $x$ and column $y$ of a Young diagram, we denote $\alpha=(x,y)$ and call $x$ and $y$ the coordinates of the node.

An $m$-partition, or a Young $m$-diagram, of size $n+1$ is an $m$-tuple of partitions such that the sum of their sizes equals $n+1$;  e.g.
the Young $3$-diagram $\left(\Box\!\Box,\Box,\Box\right)$ represents the $3$-partition $\bigl((2),(1),(1)\bigr)$ of size $4$.

We shall understand an $m$-partition as a set of $m$-nodes, where an $m$-node $\alm$ is a pair
$\{\alpha,k\}$ consisting of a node $\alpha$ and an integer $k=1,\dots,m$, indicating to which diagram in the $m$-tuple the node belongs. The integer $k$ will be called {\it position}
of the $m$-node, and we set $\pos(\alm):=k$.

For an $m$-partition $\lam$, an $m$-node $\alm$ of $\lam$ is called {\it removable} if the set of $m$-nodes obtained from $\lam$ by removing $\alm$ is still an $m$-partition. An $m$-node $\bem$ not in $\lam$ is called {\it  addable} if the set of $m$-nodes obtained from $\lam$ by adding $\bem$ is still an $m$-partition. For an $m$-partition $\lam$, we denote by
${\cal{E}}_-(\lam)$ the set of removable $m$-nodes of $\lam$ and by ${\cal{E}}_+(\lam)$ the set of addable $m$-nodes of
$\lam$.
For example, the removable/addable $m$-nodes (marked with $-/+$) for the $3$-partition $\left(\Box\!\Box,\Box,\Box\right)$
are
\[
\left(\begin{array}{l}\fbox{$\phantom{-}$}\fbox{$-$}\fbox{$+$}\\ \fbox{$+$}\end{array}
\, ,\begin{array}{l}\fbox{$-$}\fbox{$+$}\\ \fbox{$+$}\end{array}\, ,\begin{array}{l}\fbox{$-$}\fbox{$+$}\\ \fbox{$+$}\end{array}\ \right).
\]

Let $\lam$ be an $m$-diagram of size $n+1$. A standard $m$-tableau of shape $\lam$ is obtained by placing the numbers $1,\dots,n+1$ in the $m$-nodes of the diagrams of $\lam$ in such a way that the numbers in the nodes ascend along rows and down columns in every diagram.
The \emph{size} of a standard $m$-tableau is the size of its shape.

Let $q,v_1,\dots,v_m$ be the parameters of the cyclotomic Hecke algebra $H(m,1,n+1)$ and let $\alm=\{\alpha,k\}$
be an $m$-node with $\alpha=(x,y)$.
We denote by $\cc(\alm)$ the classical content of the node $\alpha$, $\cc(\alm):=y-x$, and by $c(\alm)$ the {\it quantum content} of the $m$-node $\alm$,  $c(\alm):=v_kq^{2\cc(\alm)}=v_kq^{2(y-x)}$.

For a standard $m$-tableau ${\cal{T}}$ of shape $\lam$ let $\alm_i$ be the $m$-node of ${\cal{T}}$ occupied by the number $i$, $i=1,\dots,n+1$;
we set $c({\cal{T}}|i):=c(\alm_i)$, $\cc({\cal{T}}|i):=\cc(\alm_i)$ and $\pos({\cal{T}}|i):=\pos(\alm_i)$.
For example, for the standard $3$-tableau ${\cal{T}}^{^{\phantom{A}}}\!\!\!\!={\textrm{\tiny{$\left(
\fbox{\scriptsize{$1$}}\fbox{\scriptsize{$3$}}\, ,\fbox{\scriptsize{$2$}}\, ,\fbox{\scriptsize{$4$}}\right)$}}}$ we have
\begin{alignat*}{5}
& c({\cal{T}}|1)=v_1 ,\qquad && c({\cal{T}}|2)=v_2 ,\qquad && c({\cal{T}}|3)=v_1q^2\qquad && \text{and}\qquad  c({\cal{T}}|4)=v_3 ,& \\
& \cc({\cal{T}}|1)=0 ,\qquad && \cc({\cal{T}}|2)=0 ,\qquad  && \cc({\cal{T}}|3)=1\qquad && \text{and}\qquad \cc({\cal{T}}|4)=0 ,& \\
& \pos({\cal{T}}|1)=1 ,\qquad && \pos({\cal{T}}|2)=2 ,\qquad && \pos({\cal{T}}|3)=1\qquad &&\text{and}\qquad \pos({\cal{T}}|4)=3 ,&
  \end{alignat*}

{\bf Generalized hook length.} The hook of a node $\alpha$ of a partition $\lambda$ is the set of nodes of~$\lambda$ consisting of the node $\alpha$ and the nodes which lie either under $\alpha$ in the same column or to the right of $\alpha$ in the same row; the hook length $h_{\lambda}(\alpha)$ of $\alpha$ is the cardinality of the hook of $\alpha$.
We extend this def\/inition to $m$-nodes. For an $m$-node $\alm=\{\alpha,k\}$
of an $m$-partition $\lam$, the hook length of $\alm$ in $\lam$, which we denote by $h_{\lam}(\alm)$, is the hook length of the node $\alpha$ in the $k$-th partition of $\lam$.

Let $\lam$ be an $m$-partition. For $j=1,\dots,m$, let $\mathfrak{l}_{\lam,x,j}$ be the number of nodes in the line~$x$ of the $j$-th diagram of~$\lam$, and $\mathfrak{c}_{\lam,y,j}$ be the number of nodes in the column $y$ of the $j$-th diagram of~$\lam$. The hook length of an $m$-node $\alm=\{(x,y),k\}$ of $\lam$
can be rewritten as
\[
h_{\lam}(\alm)=\mathfrak{l}_{\lam,x,k}+\mathfrak{c}_{\lam,y,k}-x-y+1
 .
 \]
Def\/ine the generalized hook length of $\alm$ (see also \cite{CJ}) by
\[
h^{(j)}_{\lam}(\alm):=\mathfrak{l}_{\lam,x,j}+\mathfrak{c}_{\lam,y,k}-x-y+1\qquad  \text{for}\  \ j=1,\dots,m  ;
\]
in particular, $h^{(k)}_{\lam}(\alm)=h_{\lam}(\alm)$ is the usual hook length.

For an $m$-partition $\lam$, we def\/ine
\begin{gather}\label{m-croc2}
\textsf{F}_{\lam}=\prod_{\alm\in \lam}\left(\frac{q^{\cc(\alm)}}{\left[h_{\lam}(\alm)\right]_q}\prod_{\textrm{\scriptsize{$\begin{array}{c}k=1,\dots,m\\k\neq\pos(\alm)\end{array}$}}}
\frac{q^{-\cc(\alm)}}{v_{\pos(\alm)}q^{-h^{(k)}_{\lam}(\alm)}-v_kq^{h^{(k)}_{\lam}(\alm)}}\right) ,
\end{gather}
where $[j]_q:=q^{j-1}+q^{j-3}+\cdots +q^{-j+1}$ for a non-negative integer $j$. Under the restrictions (\ref{sesi1})--(\ref{sesi3}),
the number $\textsf{F}_{\lam}$ is well def\/ined for any $m$-partition $\lam$ of size less or equal to $n+1$ since
$h_{\lam}(\alm)\leq n+1$ and $h^{(k)}_{\lam}(\alm)\leq n$ if $k\neq\pos(\alm)$ for any $\alm\in \lam$.

\section[Idempotents and Jucys--Murphy elements of $H(m,1,n+1)$]{Idempotents and Jucys--Murphy elements of $\boldsymbol{H(m,1,n+1)}$}\label{sec-idem}

In this section we recall the def\/inition and some properties, from  \cite{AK},
of the Jucys--Murphy elements of the algebra $H(m,1,n+1)$, together with some facts about an explicit realization of the irreducible representations of $H(m,1,n+1)$. We then derive, in the same spirit as in~\cite{Mo}, an inductive formula, that we will use in the next section, for the primitive idempotents corresponding to this realization.

The Jucys--Murphy elements $J_i$, $i=1,\dots,n+1$, of the algebra $H(m,1,n+1)$ are def\/ined by the following initial condition and recursion
\begin{gather*}
J_1=\tau\qquad\textrm{and}\qquad J_{i+1}=\sigma_iJ_i\sigma_i,\qquad i=1,\dots,n.
\end{gather*}
We recall that, under the restrictions (\ref{sesi1})--(\ref{sesi3}), the elements $J_i$, $i=1,\dots,n+1$, form a maximal commutative set
(that is, generate a maximal commutative subalgebra) of $H(m,1,n+1)$ \cite[Proposition~3.17]{AK}.
Recall also that
\begin{gather*}
J_i\sigma_k=\sigma_kJ_i\qquad\textrm{for \ \ $k\neq i-1,i$.}
\end{gather*}

The isomorphism classes of irreducible $\mathbb{C}$-representations of $H(m,1,n+1)$ are in bijection with the set of $m$-partitions of size $n+1$.
We use the labeling and the explicit realization of the irreducible representations of $H(m,1,n+1)$ given in \cite{AK}. Namely, for any $m$-partition $\lam$ of size $n+1$, the irreducible representation $V_{\lam}$ of $H(m,1,n+1)$ corresponding to $\lam$ has a basis $\{v_{{\cal{T}}}\}$ indexed by the set of standard $m$-tableaux of shape $\lam$, and is characterized (up to a diagonal change of basis) by the fact that the Jucys--Murphy elements act diagonally by
\[
J_i(v_{{\cal{T}}})=c({\cal{T}}|i) v_{{\cal{T}}}  ,\qquad  i=1,\dots,n+1.
\]
We will not need the explicit formulas for the action of the generators of $H(m,1,n+1)$ on basis elements $v_{{\cal{T}}}$.

The restriction of irreducible  representations of $H(m,1,n+1)$ to $H(m,1,n)$ is determined
by inclusion of $m$-partitions, that is, for $H(m,1,n)$-modules, we have
\begin{gather}\label{Res}
V_{\lam}\cong\bigoplus_{\mum\subset\lam ,\ \mum\ \text{of size}\ n} V_{\mum}  .
\end{gather}
Moreover, in this decomposition, $V_{\mum}$ is the space spanned by the basis vectors $v_{{\cal{T}}}$, with ${\cal{T}}$ such that the standard $m$-tableau (of size $n$) obtained by removing from ${\cal{T}}$ the $m$-node containing $n+1$  is of shape $\mum$.

For a standard $m$-tableau ${\cal{T}}$ of size $n+1$, we denote by $E_{{\cal{T}}}$ the primitive idempotent of $H(m,1,n+1)$ corresponding to $v_{{\cal{T}}}$, uniquely def\/ined by
$E_{{\cal{T}}}v_{{\cal{T'}}}=\delta_{{\cal{T}}{\cal{T'}}}v_{{\cal{T}}}$.
The results recalled above imply that $\{E_{{\cal{T}}}\}$, where ${\cal{T}}$ runs through the set of standard $m$-tableaux of size $n+1$, is a~complete set of pairwise orthogonal primitive idempotents of $H(m,1,n+1)$. Moreover, we have
by construction
\begin{gather}\label{spec-JM1}
J_i E_{{\cal{T}}}=E_{{\cal{T}}}J_i=c({\cal{T}}|i) E_{{\cal{T}}} ,\qquad i=1,\dots,n+1  .
\end{gather}
Due to the maximality of the commutative set formed by the Jucys--Murphy elements, the idempotent $E_{{\cal{T}}}$ can be expressed in terms of the elements $J_i$, $i=1,\dots,n+1$.
Let $\gam$ be the $m$-node of ${\cal{T}}$ containing the number $n+1$. As the $m$-tableau ${\cal{T}}$ is standard, the $m$-node~$\gam$ of~$\lam$ is removable. Let ${\cal{U}}$ be the standard $m$-tableau obtained from ${\cal{T}}$ by removing the $m$-node~$\gam$, and let $\mum$ be the shape of~${\cal{U}}$.
By (\ref{Res}) and (\ref{spec-JM1}), the inductive formula for $E_{{\cal{T}}}$ in terms of the Jucys--Murphy elements reads
\begin{gather*}
E_{{\cal{T}}}=E_{{\cal{U}}}\prod_{\bem\colon \begin{array}{l}\scriptstyle{\bem\in
{\cal{E}}_+(\mum)}\\[-0.5ex] \scriptstyle{\bem\neq \gam}\end{array}} \frac{J_{n+1}-c(\bem)}{c(\gam)-c(\bem)}  ,
\end{gather*}
with the initial condition: $E_{{\cal{U}}_0}=1$ for the unique $m$-tableau ${\cal{U}}_0$ of size $0$. Here $E_{{\cal{U}}}$ is considered
as an element of the algebra $H(m,1,n+1)$. Note that, due to the restrictions (\ref{sesi1})--(\ref{sesi3}), we have $c(\bem)\neq c(\gam)$ for any $\bem\in {\cal{E}}_+(\mum)$ such that $\bem\neq \gam$.

Let $\{ {\cal{T}}_1,\dots,{\cal{T}}_a\}$ be the set of pairwise dif\/ferent standard $m$-tableaux which can be obtained from ${\cal{U}}$ by adding an $m$-node with number $n+1$.
As a consequence of (\ref{Res}), we have the formula
\begin{gather}\label{somme-idem1}
E_{{\cal{U}}}=\sum_{i=1}^a E_{{\cal{T}}_i} .
\end{gather}
The element $J_{n+1}$ satisf\/ies a polynomial equation of f\/inite order so its resolvent is well def\/ined and
\[ E_{{\cal{U}}}  \frac{u-c({\cal{T}}|n+1)}{u-J_{n+1}}\]
is a rational function in an indeterminate $u$ with values in $H(m,1,n+1)$. Replacing $E_{{\cal{U}}}$ by the right-hand side of (\ref{somme-idem1}) and using (\ref{spec-JM1}), we obtain that this function
is non-singular at $u=c({\cal{T}}|n+1)$ and moreover, due to the restrictions (\ref{sesi1})--(\ref{sesi3}),
\begin{gather}\label{idem-JM-inter1} E_{{\cal{U}}}
\frac{u-c({\cal{T}}|n+1)}{u-J_{n+1}} \Bigr\rvert_{u=c({\cal{T}}|n+1)}=E_{{\cal{T}}} .
\end{gather}

\section[Fusion formula for the algebra $H(m,1,n+1)$]{Fusion formula for the algebra $\boldsymbol{H(m,1,n+1)}$}\label{sec-fus}

In this section, we prove, in Theorem \ref{prop-fus1} below, the fusion formula for the primitive idempo\-tents~$E_{{\cal{T}}}$. We use the inductive formula (\ref{idem-JM-inter1}) for $E_{{\cal{T}}}$.

Let $\phi_{k}$, for $k=1,\dots,n+1$, be the rational functions in variables $u_1,\dots,u_k$ with values in the algebra $H(m,1,n+1)$ def\/ined by
$\phi_{1}(u_1):=\bta(u_1)$ and, for $k=1,\dots,n$,
 \begin{gather*}\phi_{k+1}(u_1,\dots,u_k,u_{k+1}):=
\bsi_k(u_{k+1},u_k)\phi_{k}(u_1,\dots,u_{k-1},u_{k+1})\sigma_k^{-1}\\
\qquad {} =   \bsi_{k}(u_{k+1},u_k)\bsi_{k-1}(u_{k+1},u_{k-1})\dots \bsi_{1}(u_{k+1},u_1)\bta(u_{k+1})\sigma_1^{-1}\dots \sigma_{k-1}^{-1}\sigma_k^{-1}  . \end{gather*}
Def\/ine the following rational function $\Phi$ in variables $u_1,\dots,u_{n+1}$ with values in $H(m,1,n+1)$:
\begin{gather*}
\Phi(u_1,\dots,u_{n+1}):=
\phi_{n+1}(u_1,\dots,u_n,u_{n+1})
\phi_{n}(u_1,\dots,u_{n-1},u_{n}) \cdots \phi_{1}(u_1) .
\end{gather*}

Let $\lam$ be an $m$-partition of size $n+1$ and ${\cal{T}}$ a standard $m$-tableau of shape $\lam$. For $i=1,\dots,n+1$, we set $c_i:=c({\cal{T}}|i)$.

\begin{Theorem}\label{prop-fus1}
The idempotent $E_{{\cal{T}}}$ corresponding to the standard $m$-tableau ${\cal{T}}$ of shape $\lam$ can be obtained by the following consecutive evaluations
\begin{gather*}
E_{{\cal{T}}}=\mathsf{F}_{\lam}\Phi(u_1,\dots,u_{n+1})\bigr\rvert_{u_1=c_1}\cdots\bigr\rvert_{u_n=c_n}\bigr\rvert_{u_{n+1}=c_{n+1}}  ,
\end{gather*}
with $\mathsf{F}_{\lam}$ defined in \eqref{m-croc2}.
\end{Theorem}

We will prove the theorem in this section in several steps.

Until the end of the text, $\gam$ and $\dem$
denote the $m$-nodes of $\cal{T}$ containing the numbers $n+1$ and $n$  respectively; ${\cal{U}}$ is the standard $m$-tableau obtained from~${\cal{T}}$ by removing~$\gam$, and~$\mum$ is the shape of~${\cal{U}}$; also, ${\cal{W}}$ is the standard $m$-tableau obtained from~${\cal{U}}$ by removing the $m$-node~$\dem$ and~$\num$ is the shape of ${\cal{W}}$.

For a standard $m$-tableau ${\cal{V}}$ of size $N$, we def\/ine
the following rational function in a variab\-le~$u$ with complex values
\begin{gather}\label{def-F1}
F_{{\cal{V}}}(u):=\frac{u-c({\cal{V}}|N)}{(u-v_1)\cdots(u-v_m)}\prod_{i=1}^{N-1}\frac{\bigl(u-c({\cal{V}}|i)\bigr)^2}
{\bigl(u-q^2c({\cal{V}}|i)\bigr)\bigl(u-q^{-2}c({\cal{V}}|i)\bigr)} ;
\end{gather}
by convention,
$F_{{\cal{V}}}(u):= \frac{u-c({\cal{V}}|1)}{(u-v_1)\cdots(u-v_m)} $
for $N=1$.

\begin{Proposition} \label{prop-idem11}
We have
\begin{gather}\label{eq-idem11}
F_{{\cal{T}}}(u)\phi_{n+1}(c_1,\dots,c_n,u)E_{{\cal{U}}}=\frac{u-c_{n+1}}{u-J_{n+1}}  E_{{\cal{U}}}.
\end{gather}
\end{Proposition}

\begin{proof} We prove (\ref{eq-idem11}) by induction on $n$.
As $J_1=\tau$, we have by (\ref{bax-t2})
\[
\frac{u-c_{1}}{u-J_{1}}=\frac{u-c_{1}}{(u-v_1)\cdots(u-v_m)}\bta(u)  ,
\]
which verif\/ies the basis of induction ($n=0$).

We have: $E_{{\cal{W}}}E_{{\cal{U}}}=E_{{\cal{U}}}$ and $E_{{\cal{W}}}$ commutes with $\sigma_n$.
Rewrite the left-hand side of (\ref{eq-idem11}) as
\[
F_{{\cal{T}}}(u)\bsi_{n}(u,c_n)\cdot\phi_n(c_1,\dots,c_{n-1},u) E_{{\cal{W}}}\cdot\sigma_n^{-1} E_{{\cal{U}}} .
 \]
By the induction hypothesis we have for the left-hand side of (\ref{eq-idem11})
\[
F_{{\cal{T}}}(u)\bigl(F_{{\cal{U}}}(u)\bigr)^{-1}\bsi_{n}(u,c_n)\frac{u-c_n}{u-J_n}\sigma_n^{-1} E_{{\cal{U}}} .
\]
Since $J_{n+1}$ commutes with $E_{{\cal{U}}}$, the equality (\ref{eq-idem11}) is equivalent to
\begin{gather}
F_{{\cal{T}}}(u)\bigl(F_{{\cal{U}}}(u)\bigr)^{-1}(u-c_n)\sigma_n^{-1}(u-J_{n+1})E_{{\cal{U}}}\nonumber\\
\qquad{}=
 {\frac{(u-c_{n+1})(u-c_n)^2}{(u-q^2c_n)(u-q^{-2}c_n)}}(u-J_n)\bsi_{n}(c_n,u)E_{{\cal{U}}} \label{eq-inter11}
\end{gather}
(the inverse of $\bsi_{n}(u,c_n)$ is calculated with the help of (\ref{baxt-inv1})).
By (\ref{def-F1}),
\begin{gather*}
F_{{\cal{T}}}(u)\bigl(F_{{\cal{U}}}(u)\bigr)^{-1}(u-c_n)=(u-c_{n+1})\frac{(u-c_n)^2}{(u-q^2c_n)(u-q^{-2}c_n)}.
\end{gather*}
Therefore, to prove (\ref{eq-inter11}), it remains to show that
\begin{gather}\label{eq-inter21}
\sigma_n^{-1}(u-J_{n+1})E_{{\cal{U}}}=(u-J_n)\bsi_{n}(c_n,u)E_{{\cal{U}}}.
\end{gather}
Replacing $J_{n+1}$ by $\sigma_nJ_n\sigma_n$, we write the left-hand side of (\ref{eq-inter21}) in the form
\begin{gather}\label{eq-inter31}
\big(u\sigma_n^{-1}-J_{n}\sigma_n\big)E_{{\cal{U}}}  .
\end{gather}
As  $J_nE_{{\cal{U}}}=c_nE_{{\cal{U}}}$, the right-hand side of (\ref{eq-inter21}) is
\begin{gather*}
\left(u\sigma_n-J_n\sigma_n+\big(q-q^{-1}\big)(u-c_n)\frac{u}{c_n-u}\right)E_{{\cal{U}}}
\end{gather*}
and thus coincides with  (\ref{eq-inter31}).
\end{proof}

To prove Theorem~\ref{prop-fus1}, we need the following
information about the behavior of the rational function $F_{{\cal{T}}}(u)$ at $u=c_{n+1}$.
\begin{Proposition}\label{prop-f1}
The rational function $F_{{\cal{T}}}(u)$ is non-singular at $u=c_{n+1}$, and moreover
\begin{gather*}
F_{{\cal{T}}}(c_{n+1})=\mathsf{F}^{\phantom{-1}}_{\lam}  \mathsf{F}_{\mum}^{-1} ,
\end{gather*}
\end{Proposition}

We will prove this proposition with the help of Lemmas \ref{lemm-f1b} and \ref{lemm-f1} below, which involve the combinatorics of multi-partitions.

\begin{Lemma}\label{lemm-f1b} We have
\begin{gather}\label{eq-F1b}
F_{{\cal{T}}}(u)=(u-c_{n+1}) {\prod_{\bem\in {\cal{E}}_-(\mum)}\left(u-c(\bem)\right)} {\prod_{\alm\in {\cal{E}}_+(\mum)}\left(u-c(\alm)\right)^{-1}}  .
\end{gather}
\end{Lemma}

\begin{proof}
 The proof  is by induction on $n$. For $n=0$, we have
\[
F_{{\cal{T}}}(u)=\frac{u-c_1}{(u-v_1)\cdots(u-v_m)},
\]
which is equal to the right-hand side of (\ref{eq-F1b}).

Now, for $n>0$, we rewrite (\ref{def-F1}) for ${\cal{V}}={\cal{T}}$ as
\[
F_{{\cal{T}}}(u)=\frac{u-c_{n+1}}{(u-v_1)\cdots(u-v_m)}\frac{(u-c_n)^2}{(u-q^2c_n)(u-q^{-2}c_n)}\prod_{i=1}^{n-1}\frac{(u-c_i)^2}{(u-q^2c_i)(u-q^{-2}c_i)}  .
\]
Using the induction hypothesis, we obtain
\begin{gather}\label{eq-F1c}
F_{{\cal{T}}}(u)=\frac{(u-c_{n+1})(u-c_n)^2}{(u-q^2c_n)(u-q^{-2}c_n)}
{\prod_{\bem\in {\cal{E}}_-(\num)}\left(u-c(\bem)\right)} {\prod_{\alm\in {\cal{E}}_+(\num)}\left(u-c(\alm)\right)^{-1}}  .
\end{gather}
Denote by $\dem_t$ and $\dem_b$ the $m$-nodes which are, respectively, just above and just below~$\dem$,~$\dem_l$ and~$\dem_r$ the $m$-nodes which are, respectively, just on the left and just on the right of $\dem$; it might happen that one of the coordinates of $\dem_t$ (or $\dem_l$) is not positive, and in this situation, by def\/inition, $\dem_t\notin {\cal{E}}_-(\num)$ (or $\dem_l\notin {\cal{E}}_-(\num)$).
It is straightforward to see that:
\begin{itemize}\itemsep=0pt
\item If $\dem_t,\dem_l\notin {\cal{E}}_-(\num)$ then \[{\cal{E}}_-(\mum)={\cal{E}}_-(\num)\cup\{\dem\}\quad\textrm{and}\quad {\cal{E}}_+(\mum)=\left({\cal{E}}_+(\num)\cup\{\dem_b,\dem_r\}\right)\backslash\{\dem\}\ .\]

\item If $\dem_t\in {\cal{E}}_-(\num)$ and $\dem_l\notin {\cal{E}}_-(\num)$ then \[{\cal{E}}_-(\mum)=\left({\cal{E}}_-(\num)\cup\{\dem\}\right)\backslash\{\dem_t\}\qquad\textrm{and}\qquad {\cal{E}}_+(\mum)=\left({\cal{E}}_+(\num)\cup\{\dem_b\}\right)\backslash\{\dem\}  .\]

\item If $\dem_t\notin {\cal{E}}_-(\num)$ and $\dem_l\in {\cal{E}}_-(\num)$ then \[{\cal{E}}_-(\mum)=\left({\cal{E}}_-(\num)\cup\{\dem\}\right)\backslash\{\dem_l\}\qquad\textrm{and}\qquad {\cal{E}}_+(\mum)=\left({\cal{E}}_+(\num)\cup\{\dem_r\}\right)\backslash\{\dem\}  .\]

\item If $\dem_t,\dem_l\in {\cal{E}}_-(\num)$ then \[{\cal{E}}_-(\mum)=\left({\cal{E}}_-(\num)\cup\{\dem\}\right)\backslash\{\dem_t,\dem_l\}\qquad\textrm{and}\qquad {\cal{E}}_+(\mum)={\cal{E}}_+(\num)\backslash\{\dem\}  .\]
\end{itemize}
In each case,
using that $c(\dem_t)=c(\dem_r)=q^2c_n$ and $c(\dem_b)=c(\dem_l)=q^{-2}c_n$,
it follows that the right-hand side of (\ref{eq-F1c}) is equal to
\[(u-c_{n+1}) {\prod_{\bem\in {\cal{E}}_-(\mum)}\left(u-c(\bem)\right)} {\prod_{\alm\in {\cal{E}}_+(\mum)}\left(u-c(\alm)\right)^{-1}} ,
\]
which establishes the formula (\ref{eq-F1b}).
\end{proof}

\begin{Lemma}\label{lemm-f1}
We have
\begin{gather*}
\prod_{\bem\in {\cal{E}}_-(\mum)}\left(c_{n+1}-c(\bem)\right)\prod_{\alm\in {\cal{E}}_+(\mum)\backslash\{\gam\}}\left(c_{n+1}-c(\alm)\right)^{-1}=\mathsf{F}_{\lam} \mathsf{F}_{\mum}^{-1}  .
\end{gather*}
\end{Lemma}

\begin{proof}
{\bf 1.} The def\/inition (\ref{m-croc2}),  for a partition $\lambda$,  reduces to
\[
\textsf{F}_{\lambda}:=\prod_{\alpha\in \lambda}\frac{q^{\cc(\alpha)}}{\left[h_{\lambda}(\alpha)\right]_q} .
\]
The Lemma \ref{lemm-f1} for a partition $\lambda$ is established in \cite[Lemma 3.2]{IMO}.

{\bf 2.} Set $k=\text{pos}(\gam)$. Def\/ine, for an $m$-partition $\tem$,
\begin{gather*}
\widetilde{\textsf{F}}
_{\tem}:=\prod_{\alm\in \tem}\frac{q^{\cc(\alm)}}{\left[h_{\tem}(\alm)\right]_q}  ,
\end{gather*}
and, for $j=1,\dots,m$ such that $j\neq k$,
\begin{gather}\label{eq-F11b}
\textsf{F}^{(j)}_{\tem}:= \prod_{\textrm{\scriptsize{$\begin{array}{c}\alm\in\tem\\ \pos(\alm)=k\end{array}$}}} \frac{q^{-\cc(\alm)}}{v_kq^{-h^{(j)}_{\tem}(\alm)} -v_jq^{h^{(j)}_{\tem}(\alm)}} \prod_{\textrm{\scriptsize{$\begin{array}{c}\alm\in\tem\\ \pos(\alm)=j\end{array}$}}} \frac{q^{-\cc(\alm)}}{v_jq^{-h^{(k)}_{\tem}(\alm)} -v_kq^{h^{(k)}_{\tem}(\alm)}}  .
\end{gather}
By (\ref{m-croc2}), we have
\begin{gather}\label{eq-F11c}
\textsf{F}_{\tem}=\widetilde{\textsf{F}}_{\tem}\prod_{\textrm{\scriptsize{$\begin{array}{c}j=1,\dots,m\\j\neq k\end{array}$}}}\textsf{F}^{(j)}_{\tem}  .
\end{gather}
Fix $j\in\{1,\dots,m\}$ such that $j\neq k$. We shall show that
\begin{gather}\label{eq-F11d}
 \prod_{\textrm{\scriptsize{$\begin{array}{c}\bem\in {\cal{E}}_-(\mum)\\ \pos(\bem)=j\end{array}$}}} \left(c_{n+1}-c(\bem)\right) \prod_{\textrm{\scriptsize{$\begin{array}{c}\alm\in {\cal{E}}_+(\mum)\backslash\{\gam\}\\ \pos(\alm)=j\end{array}$}}} \left(c_{n+1}-c(\alm)\right)^{-1}=\textsf{F}^{(j)}_{\lam}\big(\textsf{F}^{(j)}_{\mum}\big)^{-1}  .
\end{gather}

Let $p_1<p_2<\cdots<p_s$ be positive integers such that the $j$-th partition of $\mum$ is $(\mu_1,\dots,\mu_{p_s})$ with
\[
\mu_1=\cdots=\mu_{p_1}>\mu_{p_1+1}=\cdots=\mu_{p_2}>\cdots >\mu_{p_{s-1}+1}=\cdots=\mu_{p_s}>0  .
\]
We set $p_0:=0$, $p_{s+1}:=+\infty$ and $\mu_{p_{s+1}}:=0$.
Assume that the $m$-node $\gam$ lies in the line $x$ and column $y$.
The left-hand side of (\ref{eq-F11d}) is equal to
\begin{gather}\label{eq-F11dg}
\prod_{b=1}^s\big(v_kq^{2(y-x)}-v_jq^{2(\mu_{p_b}-p_b)}\big)  \prod_{b=1}^{s+1}\big(v_kq^{2(y-x)}-v_jq^{2(\mu_{p_b}-p_{b-1})}\big)^{-1} .
\end{gather}

The factors in the product (\ref{eq-F11b}) correspond to
$m$-nodes of an $m$-partition.
The $m$-nodes  lying neither in the column $y$ of the $k$-th diagrams (of $\lam$ or $\mum$) nor in the line $x$ of the $j$-th
diagrams do not contribute to the right-hand side of~(\ref{eq-F11d}).
Let $t\in\{0,\dots,s\}$ be such that $p_t<x\leq p_{t+1}$. The contribution from the $m$-nodes in the column $y$ and lines $1,\dots,p_t$ of the $k$-th diagrams is
\[
\prod_{b=1}^t\left(\prod_{a=p_{b-1}+1}^{p_b}\frac{v_kq^{-(\mu_{p_b}-y+x-a)}-v_jq^{(\mu_{p_b}-y+x-a)}}
{v_kq^{-(\mu_{p_b}-y+x-a+1)}-v_jq^{(\mu_{p_b}-y+x-a+1)}}\right);
\]
the contribution from the $m$-nodes in the column $y$ and lines $p_t+1,\dots,x$ of the $k$-th diagrams is
\[
\prod_{a=p_t+1}^{x-1}\left(\frac{v_kq^{-(\mu_{p_{t+1}}-y+x-a)}-v_jq^{(\mu_{p_{t+1}}-y+x-a)}}
{v_kq^{-(\mu_{p_{t+1}}-y+x-a+1)}-v_jq^{(\mu_{p_{t+1}}-y+x-a+1)}}\right) \frac{q^{-\cc(\gam)}}{v_kq^{-(\mu_{p_{t+1}}-y+1)}-v_jq^{(\mu_{p_{t+1}}-y+1)}} .
\]
The contribution from the $m$-nodes lying in the line $x$ of the $j$-th diagrams is
\[
\prod_{b=t+1}^s\prod_{a=\mu_{p_{b+1}}+1}^{\mu_{p_b}}\frac{v_jq^{-(y-a+p_b-x)}-v_kq^{(y-a+p_b-x)}}{v_jq^{-(y-a+p_b-x+1)}-v_kq^{(y-a+p_b-x+1)}} .
\]
After straightforward simplif\/ications, we obtain for the right-hand side of (\ref{eq-F11d})
\begin{gather}
q^{x-y}\prod_{b=1}^s\bigl(v_kq^{-(\mu_{p_b}-y+x-p_b)} -v_jq^{(\mu_{p_b}-y+x-p_b)}\bigr)\nonumber\\
\qquad {}\times \prod_{b=1}^{s+1}\bigl(v_kq^{-(\mu_{p_b}-y+x-p_{b-1})} -v_jq^{(\mu_{p_b}-y+x-p_{b-1})}\bigr)^{-1} .\label{eq-F11dr}
\end{gather}
The comparison of (\ref{eq-F11dg}) and (\ref{eq-F11dr}) concludes the proof of the formula (\ref{eq-F11d}).

 {\bf 3.} The assertion of the lemma is a consequence of the formulas (\ref{eq-F11c}),
 (\ref{eq-F11d}) together with the part {\bf 1} of the proof.
\end{proof}

\begin{proof}[Proof of the Proposition~\ref{prop-f1}] The formula (\ref{eq-F1b}) shows that the rational function $F_{{\cal{T}}}(u)$ is non-singular at $u=c_{n+1}$, and moreover
\begin{gather*}
F_{{\cal{T}}}(c_{n+1})= \prod_{\bem\in {\cal{E}}_-(\mum)}\left(c_{n+1}-c(\bem)\right)\prod_{\alm\in {\cal{E}}_+(\mum)\backslash\{\gam\}}\left(c_{n+1}-c(\alm)\right)^{-1} .
\end{gather*}
We use the Lemma~\ref{lemm-f1} to conclude the proof of the proposition.
\end{proof}

\begin{proof}[Proof of Theorem  \ref{prop-fus1}]
 The theorem follows, by induction on $n$, from the formula~(\ref{idem-JM-inter1}) together with Propositions \ref{prop-idem11} and \ref{prop-f1}.
\end{proof}

\begin{Example}
Consider, for $m=2$, the standard 2-tableau ${\textrm{\tiny{$\left(
\fbox{\scriptsize{$1$}}\fbox{\scriptsize{$3$}}\, ,\fbox{\scriptsize{$2$}}\right)$}}}$. The idempotent of the algebra $H(2,1,3)$ corresponding to this standard $2$-tableau reads, by the Theorem \ref{prop-fus1},
\[
\frac{\bsi_2\big(v_1q^2,v_2\big) \bsi_1\big(v_1q^2,v_1\big) \bta\big(v_1q^2\big) \sigma_1^{-1} \sigma_2^{-1} \bsi_1(v_2,v_1) \bta(v_2)
\sigma_1^{-1} \bta(v_1)}{\big(q+q^{-1}\big)\big(v_1q^{-1}-v_2q\big)(v_1-v_2)\big(v_2q^{-2}-v_1q^2\big)}  .
\]
\end{Example}

\section{Remarks on the classical limit}\label{sec-cla-lim}

Recall that the group ring $\mathbb{C}G(m,1,n+1)$ of the complex ref\/lection group $G(m,1,n+1)$ is obtained by taking the classical limit: $q\mapsto\pm1$ and $v_i\mapsto\xi_i$, $i=1,\dots,m$, where $\{\xi_1,\dots,\xi_m\}$ is the set of distinct $m$-th roots of unity.
The ``classical limit'' of the generators $\tau, \sigma_1, \dots, \sigma_{n}$ of $H(m,1,n+1)$ we denote by $t, s_1,\dots, s_{n}$.

{\bf 1.} Consider the Baxterized elements (\ref{bax-s}) with spectral parameters of the form $v_pq^{2a}$ and $v_{p'}q^{2a'}$ with
$p,p'\in\{1,\dots,m\}$. One directly f\/inds that
\begin{gather}\label{clali-fus1}
\lim_{q\to1}\lim_{v_i\to\xi_i}\bsi_i\big(v_pq^{2a},v_{p'}q^{2a'}
\big)=s_i+\frac{\delta_{p,p'}}{a-a'}  .
\end{gather}
For the Artin generators $\tilde{s}_1,\dots,\tilde{s}_n$ of the symmetric group $S_{n+1}$, the standard Baxterized
elements are given by the rational functions
\[
\tilde{s}_i+\frac{1}{a-a'}\qquad \text{for}\ \ i=1,\dots,n .
\]
In view of (\ref{clali-fus1}), we def\/ine generalized Baxterized elements for the group $G(m,1,n+1)$ as
the following functions
\begin{gather}\label{baxt-s-m}
\overline{s}_i(p,p',a,a'):=s_i+\frac{\delta_{p,p'}}{a-a'}\qquad \text{for} \ \ i=1,\dots,n .
\end{gather}
These elements satisfy the following Yang--Baxter equation with spectral parameters
\begin{gather*}
\overline{s}_i(p,p',a,a')\overline{s}_{i+1}(p,p'',a,a'')\overline{s}_{i}(p',p'',a',a'')\\
\qquad{}=\overline{s}_{i+1}(p',p'',a',a'')
\overline{s}_{i}(p,p'',a,a'')\overline{s}_{i+1}(p,p',a,a')  .
\end{gather*}
The Baxterized elements (\ref{baxt-s-m}) have been used in \cite{OPdA} for a fusion procedure for the complex ref\/lection group $G(m,1,n+1)$.

{\bf 2.} It is immediate that
\[
\lim_{v_i\to\xi_i}\mathfrak{a}_0(r)=r^m-1\qquad \text{and}\qquad \lim_{v_i\to\xi_i}\mathfrak{a}_i(r)=r^{m-i}\qquad  \text{for \ \ $i=1,\dots,m$,}
\]
where $\mathfrak{a}_i(r)$, $i=0,\dots,m$, are def\/ined in (\ref{polalicl}). It follows from (\ref{bax-t}) that
\begin{gather}\label{clali-fus2}
\lim_{v_i\to\xi_i}\bta(r)=\sum_{i=0}^{m-1}r^{m-1-i}t^i  .
\end{gather}
The rational function $\overline{t}$ def\/ined by $\overline{t}(r):=\frac{1}{m}\sum\limits_{i=0}^{m-1}r^{m-i}t^i$ with values in $\mathbb{C}G(m,1,n+1)$ was used in~\cite{OPdA} for a fusion procedure for the complex ref\/lection group $G(m,1,n+1)$.

{\bf 3.} Def\/ine, for an $m$-partition $\lam$,
\begin{gather*}
f_{\lam}:=\left(\prod_{\alm\in\lam}h_{\lam}(\alm)\right)^{-1}.
\end{gather*}
The classical limit of $\mathsf{F}_{\lam}$ is proportional to $f_{\lam}$. More precisely, we have
\begin{gather}\label{lim-m-croc}
\lim_{q\to1}\lim_{v_i\to\xi_i}\mathsf{F}_{\lam}=\mathfrak{x}_{\lam}f_{\lam}  ,\qquad
\text{where}  \ \ \mathfrak{x}_{\lam}=\frac{1}{m^n}\prod_{\alm\in\lam} \xi_{\pos(\alm)}  .
\end{gather}
The formula (\ref{lim-m-croc}) is obtained directly from (\ref{m-croc2}) since
\[
\prod\limits_{\textrm{\scriptsize{$\begin{array}{c}i=1\\[-0.2em]i\neq k\end{array}$}}}^m (\xi_k-\xi_i)=m/\xi_k \qquad \text{for} \ \ k=1,\dots,m.
\]

{\bf 4.} Using formulas (\ref{clali-fus1}), (\ref{clali-fus2}) and (\ref{lim-m-croc}), it is straightforward to check that the classical limit of the fusion procedure for $H(m,1,n+1)$ given by the Theorem \ref{prop-fus1} leads to the fusion procedure \cite{OPdA} for the group $G(m,1,n+1)$. Also, for $m=1$,
Theorem \ref{prop-fus1} coincides with the fusion procedure \cite{IMO} for the Hecke algebra and, in the classical limit, with the fusion procedure \cite{Mo} for the symmetric group.

\subsection*{Acknowledgements}
We thank the anonymous referees for valuable suggestions.

\pdfbookmark[1]{References}{ref}
\LastPageEnding

\end{document}